%
%
%
\input amstex
\documentstyle{amsppt}
\PSAMSFonts

\pageheight{18.5cm}
\magnification=\magstep1
\frenchspacing


\loadbold

\def\phi{\varphi}

\def\nin{\newline\indent}

\def\ZZ{{\Bbb Z}}

\def\qqed{{\hfill\hfill\qed}}

\def\se#1,#2,#3;{(ab)^{#1}\, a_x^{\,#2}b_x^{\,#3}}

\def\see#1,#2,#3;{(ab)^{#3} (bc)^{#1} (ca)^{#2}\,
\ell_a^{p-#2-#3}\ell_b^{q-#1-#3}\ell_c^{r-#1-#2}}

\def\seep#1,#2,#3;{(ab)^{#3} (bc)^{#1} (ca)^{#2}\,
\ell_a^{p'-#2-#3}\ell_b^{q'-#1-#3}\ell_c^{r'-#1-#2}}

\def\seed#1,#2,#3;{(ab)^{#3} (bc)^{#1} (ca)^{#2}\,
a_x^{p_a-#2-#3}b_x^{p_b-#1-#3}c_x^{p_c-#1-#2}}

\def\nse#1,#2,#3,#4,#5,#6;{(#1 #2)(#3 #4){#5}_x{#6}_x}
\def\nsee#1,#2,#3,#4,#5,#6;{(#1 #2)(#3 #4)#5_x#6_x}
\def\nsse#1,#2,#3,#4,#5,#6,#7,#8;{(#1 #2)(#3 #4)(#5 #6)(#7 #8)}

\def\tv#1,#2,#3;{[#1,#2]_{#3}}

\TagsOnRight

\input xy
\xyoption{all}
\def\ablue{\ar@{-}}
\def\ared{\ar@{.}}
\def\agreen{\ar@{~}}
\def\ayellow{\ar@{--}}

\topmatter
\title Generic decompositions and  semi-invariants for string algebras
\endtitle
\date November 4, 2009 \enddate
\author Witold Kra\'skiewicz  and Jerzy Weyman \endauthor
\address Nicholas Copernicus University
\nin Toru\'n, Poland
\endaddress
\email wkras\@mat.umk.pl \endemail
\address Department of Mathematics, Northeastern University
\nin 360 Huntington Avenue,  BOSTON,  MA 02115, USA \endaddress
\email j.weyman\@neu.edu \endemail

\thanks The  second author was partially supported by NSF grant DMS-0901185
\endthanks
\rightheadtext{Semi-invariants of string algebras}
\abstract In this paper we investigate the rings of semi-invariants for tame string algebras $A(n)$
of non-polynomial growth. We are interested in dimension vectors of band modules.
We use geometric technique related to the description of
coordinate rings on varieties of complexes. The fascinating combinatorics emerges, showing that
our rings of invariants are the rings of some toric varieties. We show that for $n\le 6$ the rings of semi-invariants are
complete intersections but we show an example for $n=7$ that this is not the case in general.
\endabstract
\endtopmatter

\document

\head  Introduction \endhead

The rings of semi-invariants of quivers allow an explicit description only for special quivers, so-called finite type and tame quivers. For wild quivers one does not expect the explicit description to be possible for every dimension vector.

For quivers without relations the structure of rings of semi-invariants reflects the representation type of a quiver. It was shown in \cite{SW} that for a quiver $Q$ of finite representation type (a Dynkin quiver) for each dimension vector the ring $SI(Q,\alpha )$ is a polynomial ring. For a tame quiver $Q$ (i.e. an extended Dynkin quiver) the rings $SI(Q,\alpha )$ are either polynomial rings or hypersurfaces. For every wild quiver one can find a dimension vector $\alpha$ such that $SI(Q,\alpha )$ is not a complete intersection.

For quivers with relations the situation is not understood well. The basic problem that arises is that representation spaces are no longer irreducible. Thus it is natural to look at the rings of semi-invariants of irreducible components of representation spaces. In this setup the correspondance from \cite{SW} is not true.
It was shown in \cite{K} that there that for a quiver with relations of finite representation type, the ring of semi-invariants on an irreducible component of a representation space does not need to be a polynomial ring.

Still,  for finite representation type quivers with relations all the rings of semi-invariants are the semigroup rings, so they have a reasonable structure.

For tame type quivers with relations an explicit description of semi-invariants seems possible.  For a class of tame concealed and tubular algebras Kristin Webster proved in \cite{We}
that the structure of semi-invariants is similar to the semi-invariants of extended Dynkin quivers.

However the class of tame quivers with relations is much more diverse. It contains, besides  the algebras  of polynomial growth (like tubular algebras),  also the algebras that are not of polynomial growth. 

In this paper we initiate the study of rings of semi-invariants for  tame algebras of nonpolynomial growth.
The important class of algebras of nonpolynomial growth is the class of  string algebras. The indecomposable representations for string algebras were calculated by
Butler and Ringel \cite{BR}. They consist of so-called string and band modules.

Here we deal with the particular class of string algebras $A(n)$ defined in Section 2.  These algebras are convenient because the irreducible components of their representation spaces can be classified. Moreover, their representation spaces are the products of varieties of complexes, so their components are defined by rank conditions and their coordinate rings can be explicitly described. This allows the combinatorial description of semi-invariants and their structure. It turns out this structure is very different from semi-invariants of tubular algebras, and is quite fascinating in itself.
We show that for the components of band modules the rings of semi-invariants are related to the coordinate rings of some toric varieties.

We show that for $n\le 5$ the rings of semi-invariants of the components of band modules are complete intersections and give an example showing it does not happen for general $n$.

The same technique combined with the usual techniques of quiver representations allows to find the generic decompositions of band 
components for the representation spaces of $A(n)$.  This combinatorics is similar to combinatorics of semi-invariants. 

We expect that these results will have natural generalizations to general string algebras.

\bigskip

\head  \S 1. The coordinate rings of varieties of complexes \endhead

We will need some results on varieties of complexes. The basic references are \cite{DCS}, \cite{PW}, section 1.

Let $K$ be a field of characteristic 0 and let $F_0$, $F_1$,..., $F_n$ be a sequence of vector spaces of dimensions $f_0$,..., $f_n$ respectively. Consider the variety $Compl(\underline f )$ of complexes
$$0\rightarrow F_n\buildrel{d_n}\over\rightarrow F_{n-1}\buildrel{d_{n-1}}\over\rightarrow\ldots\buildrel{d_2}\over\rightarrow F_1\buildrel{d_1}\over\rightarrow F_0.$$

This variety has irreducible components corresponding to possible ranks of $d_i$'s. More precisely, consider {\it the rank sequences} $(r_n ,\ldots ,r_1 )$ satisfying the conditions $r_i +r_{i+1}\le f_i$ for $0\le i\le n$ (with the convention $r_0 = r_{n+1}=0$). We call such rank sequence $(r_n ,\ldots ,r_1 )$ $\underline f$-maximal if there is no componentwise bigger sequence satisfying the same conditions. For every $\underline f$-maximal sequence $(r_n ,\ldots ,r_1 )$ we denote
$$C(\underline f ,r_n ,\ldots ,r_1 )=\lbrace (d_n ,\ldots ,d_1 )\in Compl(\underline f)\ |\ \forall i\  rank\ (d_i )\le r_i\ \rbrace .$$
Then we have

\proclaim{Proposition 1} The irreducible components of the variety $Compl(\underline f )$ are the varieties  $C(\underline f ,r_n ,\ldots ,r_1 )$ parametrized by $\underline f$-maximal rank sequences.
\endproclaim

Notice that there is a class of sequences $\underline f$ for which the variety $Compl (\underline f )$ has a distinguished component. Fix a sequence $\underline m = (m_n ,\ldots ,m_1)$ and set $f_i :=m_i +m_{i+1}$ for $0\le i\le n$. Then there is a  maximal $\underline f$-maximal rank sequence is $(m_n ,\ldots , m_1)$, for which the component contains exact complexes. We denote such dimension vector (and the corresponding component)  by ${\underline f}[\underline m]$.

The next results we need to review is the decomposition of the coordinate rings of the representation spaces $Compl({\underline f}[\underline m])$.

Recall that for a vector space $F$ of dimension $n$ the irreducible representations of $GL(F)=GL(n)$ are parametrized by the {\it dominant integral weights} for $GL(n)$ which are the nonincreasing sequences
$$\alpha = (\alpha_1 ,\ldots ,\alpha_n )$$
of integers. We denote the corresponding irreducible representation $S_\alpha F$. This representation is (up to a power of determinant representation) {\it the Schur functor} (comp. \cite{MD}, Appendix to chapter 1, \cite{W}, chapter 2). We denote $-\alpha$ the highest weight
$$-\alpha = (-\alpha_n ,\ldots ,-\alpha_1 ).$$
We have a functorial isomorphism
$$S_\alpha F=S_{-\alpha}F^* .$$

Recall that a {\it partition} of $d$ is a sequence of nonnegative integers $\lambda = (\lambda_1 ,\ldots ,\lambda_s )$ such that $\lambda_1\ge\ldots\ge\lambda_s\ge 0$ and $d=\lambda_1 +\ldots +\lambda_s$. For a partition $\lambda =(\lambda_1 ,\ldots ,\lambda_s )$ we define the dominant integral weight $-\lambda = (-\lambda_s ,\ldots ,-\lambda_1 )$.

Sometimes we will use the Littlewood-Richardson Rule for decomposing the tensor products of the irreducible representations $S_\alpha F$. We refer to \cite{MD}, chapter 1 and \cite{W}, chapter 2 for that.

We recall some useful facts on semi-invariants. A semi-invariant {\it of weight $l$} for a rational action of $GL(F)$ is the occurrence of the Schur functor $S_{(l^n )}F$. The following lemma allows to classify semi-invariants in tensor products of Schur functors.

\proclaim {Lemma 1} Let $F=K^n$ and let $\alpha =(\alpha_1 ,\ldots ,\alpha_n )$ and $\beta =(\beta_1 ,\ldots ,\beta_n )$ be two integral dominant weights for $GL(F)$. Then the tensor product $S_\alpha F\otimes S_\beta F$ contains a semi-invariant of weight $l$ if and only if for each $i=1,\ldots ,n$ we have
$$\alpha_i +\beta_{n+1-i}=l.$$
If such case occurs, $S_\alpha F\otimes S_\beta F$ contains $S_{(l^d )}F$ with multiplicity one.
\endproclaim

Using the Schur functors one can prove

\proclaim{Proposition 2} Let us fix the rank sequence $\underline m$. 
The distinguished component  $Compl({\underline f}[\underline m])$  is irreducible and has rational singularities. 
 The coordinate ring \break $K[Compl({\underline f}[\underline m])]$ has the following decomposition
into Schur functors.
$$K[Compl({\underline f}[\underline m])]=\oplus_{\lambda_1 ,\ldots ,\lambda_n} S_{\lambda_n}F_n \otimes S_{(\lambda_{n-1},-\lambda_n )}F_{n-1}\otimes\ldots\otimes S_{(\lambda_1 ,-\lambda_2 )}F_1\otimes S_{(-\lambda_1 ) }F_0$$
where we sum over the $n$-tuples of partitions $(\lambda_n ,\ldots ,\lambda_1 )$ with
$\lambda_i$ having $m_{i-1}$ parts.
\endproclaim 

\demo{Proof} In characteristic zero this is proved in \cite{DCS}, \cite{PW}.  In characteristic $p>0$
the result is also true in the sense of good filtration. 
\qqed
\enddemo

\head  \S 2. The algebras $A(n)$ \endhead

The string algebras (i.e. monomial biserial algebras, cf. \cite{BR}) provide important examples of algebras for which the indecomposable modules are known explicitly. Among tame algebras of this type, the following algebras $A(n)$  play an important role, because they are the simplest examples of a tame algebras on nonpolynomial growth, i.e such algebras of tame representation type for which the number $N(\alpha )$ of one parameter families of indecomposable modules of dimension $\alpha$ does not grow polynomially with $\alpha$..

  Consider the quiver  $Q(n)= (Q_0 ,Q_1 )$ defined as follows
$$\matrix &\buildrel{a_n}\over\rightarrow&&\buildrel{a_{n-1}}\over
\rightarrow&\ldots&\buildrel{a_2}\over\rightarrow&&\buildrel{a_1}\over\rightarrow&\\
x_n&&x_{n-1}&&&&x_1&&x_0\\
&\buildrel{b_n}\over\rightarrow&&\buildrel{b_{n-1}}\over
\rightarrow&\ldots&\buildrel{b_2}\over\rightarrow&&\buildrel{b_1}\over\rightarrow&
\endmatrix$$

The ideal of relations $Rel(n)$ is a two sided ideal generated by  $a_i a_{i+1}=0$ and
$b_i b_{i+1}=0$ for $i=1,\ldots ,n-1$. A representation of the corresponding quiver with relations $Q(n)/Rel(n)$ is a sequence
of vector spaces $F_i = R(x_i )$ ( $i=0,\ldots ,n$ ) and the pairs of linear maps
$a_i = R(a_i ):F_i \rightarrow F_{i-1} $, $b_i =R(b_i ):F_i \rightarrow F_{i-1}$ 
for $i=1,\ldots ,n$, such that $a_i a_{i+1}=b_i b_{i+1}=0$ for $i=1,\ldots ,n-1$.

Denote $A(n):= KQ(n)/Rel(n)$ the corresponding finite dimensional algebra. Its Euler matrix is an $(n+1)\times (n+1)$ matrix
given by
$$E(n):= \left(\matrix 1&-2&2&-2&\ldots& (-1)^n 2\\
0&1&-2&2&\ldots&(-1)^{n-1} 2\\
0&0&1&-2&\ldots&(-1)^{n-2} 2\\
\ldots&\ldots&\ldots&\ldots&\ldots&\ldots\\
0&0&0&0&\ldots&1
\endmatrix\right) .$$

For the dimension vectors $\alpha := (\alpha_0 ,\ldots ,\alpha_n )$, $\beta := (\beta_0 ,\ldots ,\beta_n )$ we denote the Euler form by 
$$\langle \alpha , \beta \rangle := \alpha E(n)\beta^t .$$.

The Tits form given by the symmetrization of the Euler formcan be expressed as a square of a linear form.
$$q(\beta_0 ,\ldots ,\beta_n )= (\sum_{i=0}^n (-1)^i\beta_i )^2 .$$

We denote by $Rep(A(n),\beta )$ the representation space of representations of $A(n)$ of dimension vector $\beta$. The group $GL(\beta ):=\prod GL(\beta_i )$ acts on $Rep(A(n),\beta )$ in a natural way.

There are two reasons why the representations of $A(n)$ are interesting: it is a tame algebra of non-polynomial type, and the irreducible components of $Rep(A(n),\beta )$ can be described explicitly. Thus our interest is in describing the general modules in these components and see how the structure of these modules is reflected in the rings of semi-invariants with respect to the group $GL(\beta )$.

The variety $Rep(A(n),\beta )$ is a product of two varieties of complexes, one composed from maps $a_i$, the other from maps $b_i$. As stated in Proposition 1, the irreducible components of a variety of complexes are given by possible rank conditions $rank (a_i )=r_i$, $rank (b_i )=s_i$ which satisfy 
$$r_i +r_{i+1}\le \beta_i , s_i +s_{i+1} \le \beta_i \eqno{(*)} $$
 for $i=0,\ldots ,n$, and none of the ranks $r_i$, $s_i$ can be increased in a way that would still satisfy $(*)$.
We denote the corresponding component by \break $C(\beta; r_1   ,\ldots ,r_n ;s_1 ,\ldots ,s_n )$ or by $C(\beta ,{\underline r}, {\underline s})$ for short.  We call a component {\it a Schur component} if and only if the general module in that component is indecomposable. The component is {\it a real Schur component} if it is a Schur component and it has an open orbit.

Let us also describe the indecomposable $A(n)$-modules. These are so-called string and band modules.
A string module $S^c (w)$ can be described by a word $w$ of elements from $Q_0$, satisfying the condition
$w_{i+1}=w_i \pm 1$, and by $c\in\lbrace a,b\rbrace$. We construct the module $S^a (w)$ by assigning a copy of $K$ at a vertex $w_i$ for each $i$, and the maps $a$ act by connecting
copy $K(w_i )$ to $K(w_{i+1})$ for $i$ odd, and the maps $b$ act by connecting
copy $K(w_i )$ to $K(w_{i+1})$ for $i$ even. The module $S^b (w)$ is constructed by switching the role of $a$ and $b$.
Band modules form  1-parameter families. They correspond to cyclic words $w=(w_1 ,\ldots ,w_{2n})$, i.e. the words satisfying
$w_{i+1}=w_i \pm 1$, $w_1 =w_{2n}\pm 1$. The family $B(w, x)$ is constructed by assigning a copy of $K$ at a vertex $w_i$ for each $i$, and the maps $a$ act by connecting
copy $K(w_i )$ to $K(w_{i+1})$ for $i$ odd, and the maps $b$ act by connecting
copy $K(w_i )$ to $K(w_{i+1})$ for $i$ even, as well as connecting $K(w_{2n})$ to $K(w_1 )$ by the map $b$, but with a scalar $x$ attached.
Notice that we could define modules $B^a (w,x)$ and $B^b (w,x)$ but they would give the same families.

It is well known (comp. \cite{BR}) that the string and band modules give the complete description of the indecomposable modules for the algebra $A(n)$.

We close this section with some terminology. A component $C(\beta ,{\underline r}, {\underline s})$ is a {\it Schur component} if a general module $M$ in  $C(\beta ,{\underline r}, {\underline s})$ is indecomposable. A Schur component $C(\beta ,{\underline r}, {\underline s})$ is {\it real} if it contains an open orbit. Otherwise a Schur component $C(\beta ,{\underline r}, {\underline s})$ is {\it an imaginary Schur component}.

\head  \S 3. The generic decomposition of band modules for the algebras $A(n)$ \endhead

In this section we deal with the generic decompositions of band modules for the general algebra $A(n)$. This means we will deal with the dimension vectors  
$$\beta (m_{n} ,\ldots ,m_{1}):=[m_{n} ,\ldots ,m_1] = (m_{n} ,m_{n}+m_{n-1} ,\ldots m_{2}+m_{1} , m_1 ).$$

For each such dimension vector $\beta (m_{n} ,\ldots ,m_{1}):=[m_{n} ,\ldots ,m_1]$ the representation space $Rep(A(n),  \beta (m_{n} ,\ldots ,m_{1}))$ has a distinguished component \break
$C(\beta (m_{n},\ldots ,m_1); (m);(m))$. We will refer to these components as the {\it band module components} because these are the only irreducible components of representation spaces where the generic module is a direct sum of band modules.

Next we define the  up-and-down graphs and modules for the algebra $A(n)$. We will consider the double colored graphs with red and blue edges.

\proclaim {Definition 1(of up-and-down graphs)} Let $\beta, {\underline r}, {\underline s}$ be fixed. The up-and-down graph $\Gamma(\beta ; {\underline r}, {\underline s} )$ is defined as follows. The vertices of $\Gamma (\beta ; {\underline r}, {\underline s} )$ have double indexing, they are denoted
 $\lbrace e^{(i)}_1 ,\ldots ,e^{(i)}_{\beta_i}\rbrace$  for $0\le i\le n$.
The red edges connect vertices $ e^{(i)}_j $ and $e^{(i-1)}_j $ for $1\le j\le r_i$ and $n-i$ even. The blue edges connect the vertices $ e^{(i)}_{\beta_i  +1-j}$ and $ e^{(i-1)}_{\beta_{i-1}+1-j}$
for $1\le j\le s_i$ and $n-i$ even. Similarly, for $n-i$ the red edges connect the vertices $e^{(i)}_{\beta_i +1-j}$ and $ e^{(i-1)}_{\beta_{i-1}+1-j}$ for $1\le j\le r_i$ and the blue edges connect the vertices $e^{(i)}_j $ and $e^{(i-1)}_j$ for $1\le j\le s_i$.  
If $\beta$ is a band module component, we denote the corresponding up-and-down graph just by $\Gamma (\beta )$.
\endproclaim

\proclaim {Example 1} Here is the up-and-down graph for the dimension vector $(2,5,4,1)$.

 $$
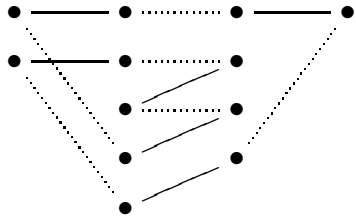
\xymatrix @R=5pt{
\bullet \ablue[r]\ared[3,1] & \bullet \ared[r] & \bullet \ablue[r] & \bullet \\
\bullet \ablue[r]\ared[3,1] & \bullet \ared[r] & \bullet \\
                   & \bullet \ared[r]\ablue[ru] & \bullet \\
                   & \bullet \ablue[ru]& \bullet\ared[-3,1] \\
                   & \bullet \ablue[ru]\\
 }$$
 \endproclaim

\proclaim {Remark 1} It is essential to observe that the connected components of an up-and-down graph are also up-and-down graphs. Also the choice of the red edges to give the upper part for even $i$ does not matter, because writing the graph upside down would make blue arrows giving the upper part for even $i$.
\endproclaim

\proclaim {Example 2} Consider the dimension vector $(5,9,7,3)$. The corresponding up-and-down graph $\Gamma (5,9,7,3)$ has two connected 
components that are isomorphic to the up-and-down graphs  $\Gamma (2,4,3,1)$ and  $\Gamma (3,5,4,2)$. Their vertices are closed and open nodes respectively
$$
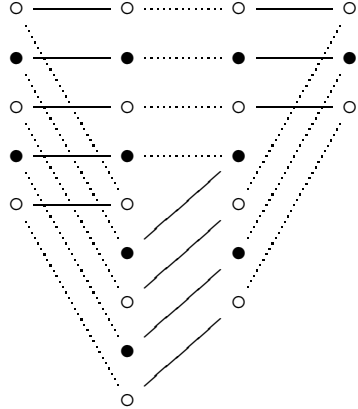
\xymatrix @R=5pt{
\circ \ablue[r]\ared[4,1] & \circ \ared[r] & \circ \ablue[r] & \circ \\
\bullet \ablue[r]\ared[4,1] & \bullet \ared[r] & \bullet \ablue[r] & \bullet \\
\circ \ablue[r]\ared[4,1] & \circ \ared[r] & \circ \ablue[r] & \circ \\
\bullet \ablue[r]\ared[4,1] & \bullet \ared[r] & \bullet \\
\circ \ablue[r]\ared[4,1] &\circ  &\circ \ared[-4,1]  &  \\
                   & \bullet \ablue[-2,1] & \bullet \ared[-4,1] \\
                   & \circ \ablue[-2,1]& \circ \ared[-4,1] \\
                   & \bullet \ablue[-2,1] & \\
                   &\circ \ablue[-2,1] &
 }$$
 \endproclaim

We define {\it a string graph} to be a graph of type $A_m$ with the edges colored in alternative way (i.e. neiboring edges have different colors). Similarly, {\it a band graph} is a graph of type ${\hat A}_m$ with the edges colored in alternative way.

\proclaim {Proposition 3}
\item{a)} For every component $C (\beta ; {\underline r}, {\underline s} )$ the graph $\Gamma (\beta ; {\underline r}, {\underline s} )$ is a direct sum of string and band graphs.
\item{b)} For the dimension vectors $\beta:=[m_{n},\ldots ,m_1 ]$ the graph $\Gamma (\beta ; {\underline r}, {\underline s} )$ is a sum of band graphs.
\endproclaim

\demo{Proof} Every vertex of $\Gamma (\beta ; {\underline r}, {\underline s} )$ is connected to at most one red and at most one blue edge. Thus the graph deomposes to strings and bands. In the case of $\beta:=[m_{n},\ldots ,m_1 ]$ every vertex is connected to exactly one red and exactly one blue edge. Thus the graph decomposes to bands.
\qqed
\enddemo

\proclaim {Remark 2} We can draw the graphs $\Gamma (\beta ; {\underline r}, {\underline s} )$ in the plane (with edges intersecting) as follows. We call the vertices $e^{(i)}_j$ the vertices on the level $i$. We draw the vertices of $\Gamma (\beta ; {\underline r}, {\underline s} )$ in the lattice $\ZZ^2$  by drawing vertex $e^{(i)}_j$ as a point with coordinates $(i,-j)$.
\endproclaim

\proclaim {Example 3} Here is a graph $\Gamma ((2,3,2,4,2);2,0,2,2;2,1,1,2)$.
$$
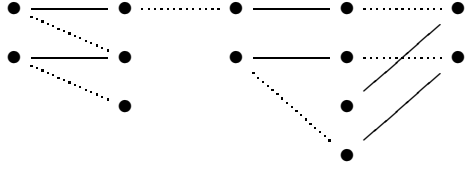
\xymatrix @R=5pt{
\bullet \ablue[r] \ared[1,1] &\bullet \ared[r] &\bullet \ablue[r]&\bullet \ared[r]&\bullet \\
\bullet \ablue[r] \ared[1,1]  &\bullet & \bullet \ablue[r]\ared[2,1] &\bullet \ared[r]&\bullet \\
&\bullet &&\bullet \ablue[-2,1] &\\
&&&\bullet \ablue[-2,1]&
}$$
There are two connected components, one is a band graph and the other is a string graph.
\endproclaim

Now we can define the families of up-and-down modules.

\proclaim {Definition 2(of up-and-down modules)} Let $\beta ; {\underline r}, {\underline s} $ be fixed. The family up-and-down modules $V:=M(\beta ; {\underline r}, {\underline s} )$ is defined as follows.
We fix bases $\lbrace e^{(i)}_1 ,\ldots ,e^{(i)}_{\beta_i}\rbrace$,   of $V(i)$ for $0\le i\le n$.
The map $a_i$ is defined as an "upper identity" i.e. $a_i (e^{(i)}_j )=e^{(i-1)}_j $ for $1\le j\le r_i$ and $n-i$ even. The map $b_i$ is defined as a "lower identity" $b_i (e^{(i)}_{\beta_i  +1-j})= e^{(i-1)}_{\beta_{i-1}+1-j}$
for $1\le j\le s_i$ and $n-i$ even. Similarly, for $n-i$ odd $a_i$ is the "lower identity" $a_i (e^{(i)}_{\beta_i +1-j})= e^{(i-1)}_{\beta_{i-1}+1-j}$ for $1\le j\le r_i$ and $b_i$ is the "upper identity"
$b_i (e^{(i)}_j )=e^{(i-1)}_j$ for $1\le j\le s_i$.  We modify this definition by adding a scalar (in the upper most left place) in every component of $\Gamma (\beta ; {\underline r}, {\underline s} )$ that is a band.
The family of modules defined in such way is indeed a family of $A(n)$-modules. We refer to them as an "up-and-down" modules corresponding to the component 
 $C (\beta ; {\underline r}, {\underline s} )$. 
 \endproclaim

\proclaim {Remark 3} The name "up-and-down" refers to the fact that in the graphs $\Gamma (\beta ; {\underline r}, {\underline s} )$ the red (resp. blue) edges connect either the $r_i$ (resp. $s_i$) top vertices from the $i$-th and $i-1$-st set of the bottom $r_i$ (resp. $s_i$) ones.
\endproclaim

Here are the results relating up-and-down modules to generic decompositions.

\proclaim {Theorem A} The up-and-down band modules are generic in their components, i.e. the union of their orbits contains a Zariski open set in the corresponding component.
\endproclaim

\demo{Proof} We will show that up-and-down band modules do not have selfextensions, i.e. that $Ext^i_{A(n)} (M(\beta ), M'(\beta ))=0$ for generic choices of the up-and-down modules in dimension $\beta :=[m_{n},\ldots ,m_1]$, and for all $i\ge 0$. Then the theorem follows from the result of Crawley-Boevey and Schroer \cite{CBS}.

We start with the definitions regarding graphs $\Gamma (\beta ; {\underline r}, {\underline s} )$. We call a path in this graph a upper path (resp. a lower path) if it stays in the upper (lower) part of the graph. The vertex $e^{i}_j$ is called the left extremal vertex (resp. right extremal vertex) if both edges connecting to it connect to the level $i-1$ (resp. connect to the level $i+1$).
We denote $LEXT(\Gamma (\beta ;{\underline r};{\underline s}))$ (resp. $REXT(\Gamma (\beta ;{\underline r};{\underline s}))$) the set of left (resp. right extremal boxes in $\Gamma (\beta ;{\underline r};{\underline s})$.
Now we can describe the minimal projective resolution of the up-and-down module $M(\beta )$.

\proclaim {Proposition 4} The up-and-down band module $M(\beta )$ has projective dimension $1$, and its projective resolution is
$$0\rightarrow P^{(1)}\rightarrow P^{(0)}\rightarrow M(\beta )\rightarrow 0$$
where
$$P^{(1)}=\oplus_{e^{(k)}_l\in REXT(\Gamma (\beta ;{\underline r};{\underline s}))}P_k ,$$
$$P^{(0)} =\oplus_{e^{(i)}_j\in LEXT(\Gamma (\beta ;{\underline r};{\underline s}))}P_i .$$
The matrix entry of the differential from the summand $P_k$ corresponding to $e^{(k)}_l$ to the summand $P_i$ corresponding to $e^{(i)}_j$ is nonzero if the path from $e^{(i)}_j$ to $e^{(k)}_l$ goes directly to the right. The nonzero matrix entry is given by this path, multiplied by $-1$ if it is the lower path, and multiplied by the scalar $\lambda$ if the edge marked by $\lambda$ is on the path  from $e^{(i)}_j$ to $e^{(k)}_l$. 
\endproclaim

\demo{Proof} This is clear by definitions.
\qqed
\enddemo

Now we can finish the proof of the theorem A. Consider the complex
$$0\rightarrow Hom_{A(n)} (P^{(0)} , M(\beta ))\rightarrow Hom_{A(n)} (P^{(1)} , M(\beta )) \rightarrow 0.$$
We want to show that it is an isomorphism of vector spaces. 

The  vector space $Hom_{A(n)} (P^{(0)} , M(\beta ))$ has a basis given by $e^{(i)}_u\otimes e^{(i)}_j$ where $e^{(i)}_j\in LEXT(\Gamma (\beta ;{\underline r};{\underline s}))$. The  vector space $Hom_{A(n)} (P^{(1)} , M(\beta ))$ has a basis given by $e^{(k)}_v\otimes e^{(k)}_l$ where $e^{(k)}_l\in REXT(\Gamma (\beta ;{\underline r};{\underline s}))$. We construct a graph $EXT(M(\beta ))$ with the vertices of both spaces by connecting $e^{(i)}_u\otimes e^{(i)}_j$ to $e^{(k)}_v\otimes e^{(k)}_l$ if the matrix entry from the projectives cooresponding to the extremal vertices is nonzero, and the path given by matrix entry takes $e^{(i)}_u$ to $e^{(k)}_v$. 
Notice that the edges in the graph  $EXT(M(\beta ))$ consist of either upper or lower path in $\Gamma (\beta ; {\underline r}, {\underline s} )$. We refer to them as upper and lower edges in
$EXT(M(\beta ))$. The upper and lower edges in $EXT(M(\beta ))$ alternate, i.e. if a vertex has two edges adjacent to it, one is lower and the other is upper.
The main observation is that the graph $EXT(M(\beta ))$ has band and string components. The band components consist of all vertices $e^{(i)}_j\otimes e^{(i)}_j$ and $e^{(k)}_l\otimes e^{(k)}_l$. This part of the complex is exact after choosing genric scalars on each band. The remaining components are the string components. We have to show that these string components are even, so they will give exact complexes of type $K^d\rightarrow K^d$. Arrange basis vectors in $Hom_{A(n)} (P^{(i)} , M(\beta ))$ ordering the blocks in the way the extremal vertices are ordered in the band module, and within each block as we draw 
$\Gamma (\beta ;{\underline r};{\underline s})$. Now each block of boxes has one element in a band component.
Now we make the following observation about the string components: either all their elements  occur above or below the band orbit element in their block.
The key point is that  if the orbit occurs above (resp. below) the band orbit in each block, then  its endpoints will be adjacent to the  upper (resp. lower) edge in $EXT(M(\beta ))$. Thus each string component has odd number of edges and thus even number of vertices. This proves our complex has zero homology, so for generic band modules all their Ext groups are zero.
\qqed
\enddemo

\proclaim {Example 4} We show the proof of Theorem A in graphical form on an example.
Consider the dimension vector $(2,4,5,3)$. The corresponding up-and-down graph is
$$
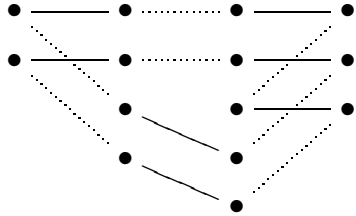
\xymatrix @R=5pt{
\bullet \ablue[r] \ared[2,1] &\bullet \ared[r] &\bullet \ablue[r] &\bullet \\
\bullet \ablue[r] \ared[2,1] &\bullet \ared[r] &\bullet\ablue[r] &\bullet \\
            &\bullet \ablue[1,1] &\bullet \ablue[r]\ared [-2,1] &\bullet \\
            &\bullet \ablue[1,1] &\bullet \ared[-2,1]&  \\
&&\bullet \ared[-2,1] &
}$$
We have the projective resolution
$$0\rightarrow P_0\oplus P_0\oplus P_0\rightarrow P_3\oplus P_1\oplus P_3\rightarrow M(\beta )_\lambda\rightarrow 0$$
where the map is given by the matrix
$$\left(\matrix \lambda a_1b_2a_3&0&b_1a_2b_3\\ -b_1&a_1&0\\ 0&-b_1a_2b_3&-a_1b_2a_3
\endmatrix
\right) .$$
The graph $EXT(M(\beta ))$ is as follows
$$
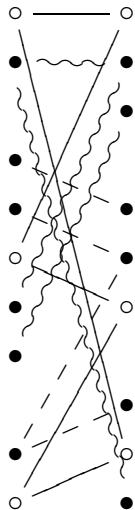
\xymatrix @R=5pt{
\circ \ablue[r] \ablue[9,1] &\circ \\
\bullet \agreen[r] \agreen[9,1]&\bullet \\
            &\bullet \\
\bullet \ayellow[1,1] &            \\
\bullet \ayellow[1,1] &\bullet \\
\circ \ablue[1,1]\ablue[-5,1] &\bullet \\
\bullet \agreen[-5,1] &\circ \\
\bullet \agreen[-5,1] &            \\
            &\bullet \\
\bullet\ayellow[-4,1]\ayellow[-1,1] &\circ \\
\circ \ablue[-4,1]\ablue[-1,1] &\bullet 
}$$
with the blocks coming from $P_3, P_1, P_3$ in the left column and the blocks coming from $P_0 ,P_0 ,P_0$ in the right column.
The band component has the vertices given as open nodes. and the edges given as 
$$\xymatrix @R=5pt{\bullet \ablue[r] &\bullet}.$$
The  components occuring above the band component have edges given as 
$$\xymatrix @R=5pt{\bullet\ayellow[r] &\bullet}.$$
The  components occurring below the band component have the edges given as
$$\xymatrix @R=5pt{\bullet\agreen[r] &\bullet}.$$
\endproclaim

Let us recall that {\it the generic decomposition} of the band module component $\alpha$ 
is a sum of band module components $\alpha^{(1)},\ldots , \alpha^{(r)}$ if a general module of dimension $\alpha$ decomposes to modules of dimensions $\alpha^{(1)},\ldots , \alpha^{(r)}$.  

We will write the generic decomposition as the decomposition of the dimension vectors. It will be understood as the decomposition of the general module into string and band modules.
In the case of the band module components one can give simple rules to calculate such decompositions. For $\alpha =[m_{n},\ldots ,m_1 ]$ we will denote by $N([m_{n},\ldots ,m_1 ])$ the number of summands in the generic decomposition of $[m_{n},\ldots ,m_1 ]$. Thus a band module component is a Schur component if and only if $N([m_{n},\ldots ,m_1 ])=1$.

\proclaim {Theorem B} The decomposition of band module components $[m_{n},\ldots ,m_1]$ is given as follows. 
\item {a)} If $m_{i+1} =m_i$ then the decomposition of $[m_{n},\ldots ,m_1]$ is the same as the decomposition of the shorter vector $[m_{n},\ldots ,m_{i+1}, m_{i-1},\ldots ,m_1]$ with the repetition of the $i$'th and $i+1$-st place added.
\item{b)} Assume that $m_i >max (m_{i-1}, m_{i+1})$. Then the decomposition of $[m_{n},\ldots ,m_1]$ is the same as that of $[m_{n},\ldots ,m_{i+1}, m_i -|m_{i-1}-m_{i+1}|, m_{i-1},\ldots ,m_1]$ meaning that the number of summands in both decompositions is the same and the summands differ only in the $i$-th place.
\item{c)} Assume that $m_i >m_{i-1}=m_{i+1}$. Then in the decomposition of $[m_{n},\ldots ,m_1]$ there are $m_i -m_{i-1}$ summands $[0,\ldots ,0,1,0,\ldots ,0]$ (with $1$ in the $i$-th place) and the remaining summands are the same as the summands in the decomposition of $[m_{n},\ldots ,m_{i+1}, m_{i+1}, m_{i-1},\ldots ,m_1]$.
\item{d)} For $n=2$ we know that if $gcd(m_2 ,m_1 )=d$, $m_i =d n_i$ then
$$[m_2 ,m_1 ]=d[n_2 ,n_1]$$
and $[n_2 ,n_1]$ is the Schur component. Thus we have $N([m_2 ,m_1 )=d$.
\item{e)} $N([m_3 ,m_2 ,m_1])= gcd (m_2 ,|m_3 -m_1 |)$. 
\endproclaim

\proclaim {Remark 4} The rules a)-d) allow to determine quickly the number of summands in the generic decomposition of band module components, and which band components are the Schur components. The slight difficulty in determining the decomposition is that in the reduction b) it is not exactly clear for which $|m_{i-1}-m_{i+1}|$ bands the coordinate $m_i$ changes. 
\endproclaim

\demo{Proof of Theorem B}. We know by Theorem A that the up-and-down band modules give generic decomposition of the band components.
 We will work with graphs $\Gamma (\beta ;{\underline r}, {\underline s})$ and will prove that their decompositions to band modules obey rules a)-d).
 Let us look at rule a). We want to show that If $m_{i+1} =m_i$ then the decomposition of $[m_{n},\ldots ,m_1]$ is the same as the decomposition of the shorter vector $[m_{n-1},\ldots ,m_{i+1}, m_{i-1},\ldots ,m_0]$ with the repetition of the $i$'th and $i+1$-st place added. Let us compare the graphs $\Gamma$ for those dimension vectors. We call them for short $\Gamma([m_0,\ldots ,m_{n-1}])$ because $\underline r =\underline s =(m_0,\ldots ,m_{n-1})$.

Let us call a vertex $u$ of $\Gamma([m_n ,\ldots ,m_{1}])$ left extremal if both its edges go to the right, and right extremal if both its edges go to the left.
Notice that the graph $\Gamma (\beta [m_{n},\ldots ,m_1])$ has no left or right extremal vertices at level $i+1$. We can now obtain the graph $\Gamma([m_n ,\ldots ,m_{i-1}, m_{i+1},\ldots ,m_{1}])$ by connecting the edges from level $i$ to level $i-1$ and following edges from level $i+1$ to level $i$  into the edges from level $i+1$ to level $i-1$, coloring them according to their left parts, and changing the coloring of the edges to the right of level $i$. Now its clear that both graphs have the same number of bands.

To prove rule b) let us assume that $m_i >max (m_{i-1}, m_{i+1})$. Compare the graphs $\Gamma([m_n,\ldots ,m_{1}])$ and $\Gamma ([m_{n},\ldots ,m_{i+1}, m_i -|m_{i-1}-m_{i+1}|, m_{i-1},\ldots ,m_1])$. They differ only between the levels $i$ and $i-1$. Assume to fix the attention that $i$ is even. The red edges between level $i+1$ and $i$ in both graphs are horizontal, while the blue edges go up  by $m_{i-1}-m_{i+1}$ (which means going down if $m_{i-1} <m_{i+1}$. The graph $\Gamma ([m_{n},\ldots ,m_{i+1}, m_i -|m_{i-1}-m_{i+1}|, m_{i-1},\ldots ,m_1])$ has $|m_{i-1}-m_{i+1}|$ vertices less on levels $i$ and $i-1$. Thus it is clear that the number of bands in both graphs is the same, with $|m_{i-1}-m_{i+1}|$ bands in $\Gamma (\beta [m_{n},\ldots ,m_1])$ having one more  zig-zag connecting $e^{(i)}_x$ to the vertex on level $i+1$ and then to the vertex $e^{(i)}_{x+m_{i+1}-m_{i-1}}$. This proves rule b). 

Finally, to prove c), assume that $m_i >m_{i-1}=m_{i+1}$. Notice that in such case both graphs differ only at levels $i+1$ and $i$. In both graphs the numbers of vertices at levels $i$ and $i-1$ are the same. In both graphs the red and blue edges between levels $i$ and $i-1$ are horizontal. But $m_i -m_{i+1}$ pairs vertices in $\Gamma([m_n,\ldots ,m_{1}])$ at levels $i$ and $i-1$ are connected by both the blue and red edge, thus leading to summands of type $[0,\ldots ,0,1,0,\ldots ,0]$ with $1$ in the $i$-th place. The rest of $\Gamma([m_n,\ldots ,m_{1}])$ is the same as $\Gamma ([m_{n-1},\ldots ,m_{i+1}, m_{i+1}, m_{i-1},\ldots ,m_0])$. This proves rule c). The proof of rule d) is straightforward, so we leave it to the reader. 
\qqed
\enddemo

We illustrate the proof of the reductions by comparing up-and-down graphs.

\proclaim {Example 5(rule a))} Compare the up-and-down graphs  $\Gamma (1,3,4,2)$ and $\Gamma (1,3,2)$.
They are
$$
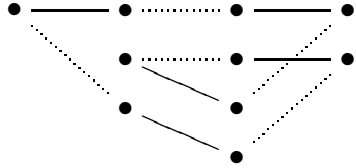
\xymatrix @R=5pt{
\bullet \ablue[r]\ared[2,1] &\bullet \ared[r] &\bullet \ablue[r]&\bullet \\
&\bullet \ared[r]\ablue[1,1]&\bullet \ablue[r]&\bullet \\
&\bullet \ablue[1,1]&\bullet \ared[-2,1]&\\
&&\bullet \ared[-2,1]&
}$$
and 
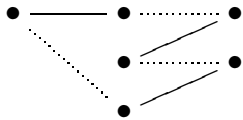
$$\xymatrix @R=5pt{
\bullet \ablue[r] \ared[2,1] &\bullet \ared[r]&\bullet \\
&\bullet \ared[r]\ablue[-1,1]&\bullet \\
&\bullet \ablue[-1,1]&
}$$
As remarked in the proof the second graph is obtained from the first by omitting the vertices coresponding to vertex 2 and replacing the paths of length 2 through these vertices by paths of length one, appropriately colored.
\endproclaim

\proclaim {Example 6(rule b))} Let us compare the up-and-down graphs  $\Gamma (1,3,9,11,4)$ and $\Gamma (1,3,7,9,4)$.
They are
$$
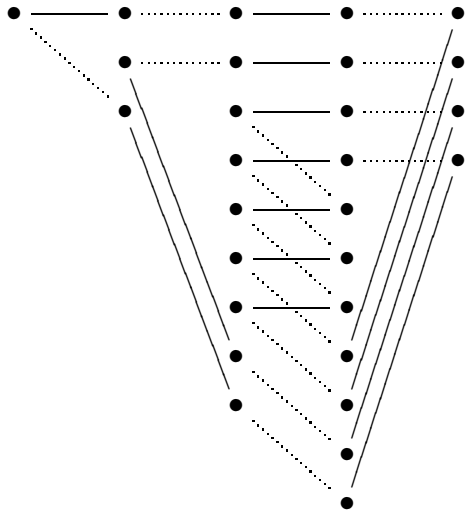
\xymatrix @R=5pt{
\bullet\ablue[r]\ared[2,1] &\bullet\ared[r] &\bullet\ablue[r] &\bullet\ared[r] &\bullet \\
&\bullet \ared[r]\ablue[6,1] &\bullet\ablue[r] &\bullet\ared[r] &\bullet \\
&\bullet \ablue[6,1] &\bullet\ablue[r]\ared[2,1] &\bullet\ared[r] &\bullet \\
&&\bullet\ablue[r]\ared[2,1] &\bullet\ared[r] &\bullet \\
&&\bullet\ablue[r]\ared[2,1] &\bullet &\\
&&\bullet\ablue[r]\ared[2,1] &\bullet &\\
&&\bullet\ablue[r]\ared[2,1] &\bullet &\\
&&\bullet\ared[2,1] &\bullet \ablue[-7,1] &\\
&&\bullet\ared[2,1] &\bullet\ablue[-7,1]  &\\
&&&\bullet\ablue[-7,1] &\\
&&&\bullet\ablue[-7,1] &
}$$
and 
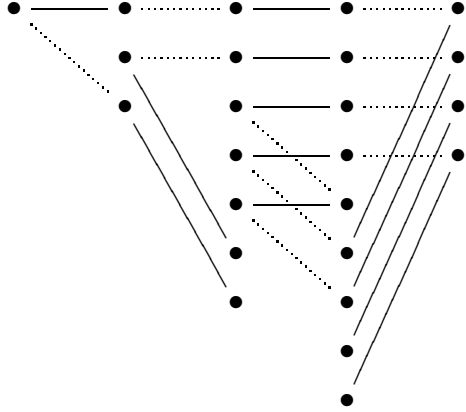
$$\xymatrix @R=5pt{
\bullet\ablue[r]\ared[2,1] &\bullet\ared[r] &\bullet\ablue[r] &\bullet\ared[r] &\bullet \\
&\bullet \ared[r]\ablue[4,1] &\bullet\ablue[r] &\bullet\ared[r] &\bullet \\
&\bullet \ablue[4,1] &\bullet\ablue[r]\ared[2,1] &\bullet\ared[r] &\bullet \\
&&\bullet\ablue[r]\ared[2,1] &\bullet\ared[r] &\bullet \\
&&\bullet\ablue[r]\ared[2,1] &\bullet &\\
&&\bullet&\bullet\ablue[-5,1] &\\
&&\bullet &\bullet\ablue[-5,1] &\\
&&&\bullet\ablue[-5,1] &\\
&&&\bullet\ablue[-5,1] &
}$$
Observe, that (as stated in the proof) both graphs have the same number (two) of connected compoenents and these components differ only 
in vertices 2 and 3 (with the components in the first graph having one more zig-zag between vertices 2 and 3). The key point in the proof is that in both graphs the edges of both colors have the same effect on the level of the vertex (the red edges are horizontal and the blue edges go down two levels).
\endproclaim

\proclaim {Example 7(rule c))} Look at the graph $\Gamma (1,4,4,1)$.
$$\xymatrix @R=5pt{
\bullet\ablue[r]\ared[3,1]&\bullet\ared[r] &\bullet\ablue[r] &\bullet \\
&\bullet\ablue[r]\ared[r] &\bullet &\\
&\bullet\ablue[r]\ared[r] &\bullet &\\
&\bullet\ablue[r] &\bullet\ared[-3,1]&
}$$
Note that in two connected components in the middle the vertices are connected by arrows of both colors. This graph differs from $\Gamma (1,2,2,1)$ by two components in the middle.
\endproclaim

The reduction rules also apply to the string components . 

\proclaim {Definition 3} A  string component is a component
$C(\beta ;{\underline r}, {\underline s})$ satisfying the condition $\sum_{i=0}^n (-1)^i\beta_i =\pm 1$, for all $i$ except one we have $r_i +r_{i-1}=\beta_i$ and for all $i$ except one we have
$s_i +s_{i-1}=\beta_i$. This means that the complex we would get from considering only maps $a_i$ would have only one dimensional homology in one spot, and similarly for a complex we would get from considering only the maps $b_i$.
\endproclaim

The reason for considering string components is that a real Schur component has to be a string component. In fact a dimension component is a string component if and only if it contains a string module.

\proclaim {Proposition 5} Let $C(\beta ;{\underline r}, {\underline s})$ be a string component. Then the graph $\Gamma(\beta ;{\underline r}, {\underline s})$ decomposes to one string graph and several band graphs. 
\endproclaim

\demo{Proof} Each vertex except two has a red and blue edge connecting to it.
\qqed
\enddemo

\proclaim {Theorem A'} Let $C(\beta ,\underline r,\underline s)$ be a string component. The corresponding up-and-down modules are generic in their components, i.e. the union of their orbits contains a Zariski open set in the corresponding component.
\endproclaim

\demo{Proof} We will show that up-and-down string modules do not have selfextensions, i.e. that $Ext^i_{A(n)} (M(\beta ), M'(\beta ))=0$ for generic choices of the up-and-down modules in dimension $\beta :=[m_{n},\ldots ,m_1]$, and for all $i\ge 0$.
Then the theorem follows form the result of Crawley-Boevey and Schroer \cite{CBS}.

We start with the definitions regarding graphs $\Gamma (\beta ; {\underline r}, {\underline s} )$. We call a path in this graph a upper path (resp. a lower path) if it stays in the upper (lower) part of the graph. The vertex $e^{i}_j$ is called the left extremal vertex (resp. right extremal vertex) if both edges connecting to it connect to the level $i+1$ (resp. connect to the level $i-1$).
We denote $LEXT(\Gamma (\beta ;{\underline r};{\underline s}))$ (resp. $REXT(\Gamma (\beta ;{\underline r};{\underline s}))$) the set of left (resp. right extremal boxes in $\Gamma (\beta ;{\underline r};{\underline s})$. The left and right end vertices of $\Gamma (\beta ; {\underline r}, {\underline s} )$ are denoted by $x_0$ and $x_\infty$. 
Now we can describe the minimal projective resolution of the up-and-down string module $M(\beta ;{\underline r};{\underline s})$.

\proclaim {Proposition 6} The up-and-down band module $M(\beta  ;{\underline r};{\underline s})$ has the projective  resolution defined as follows
$$0\rightarrow P^{(1)}\rightarrow P^{(0)}\rightarrow M(\beta  ;{\underline r};{\underline s}))\rightarrow 0$$
where
$$P^{(1)}=\oplus_{e^{(k)}_l\in REXT(\Gamma (\beta ;{\underline r};{\underline s})}P_k \oplus P^{(1,0)}\oplus P^{(1,\infty )}.,$$
$$P^{(0)} =\oplus_{e^{(i)}_j\in LEXT(\Gamma (\beta ;{\underline r};{\underline s})}P_i \oplus P^{(0,0)}\oplus P^{(0,\infty )}.$$
$$P^{(i)}= P^{(i,0)}\oplus P^{(i,\infty )}.$$
Here the terms $P^{(i,0)}$, $ P^{(i,\infty )}$ corresponding to the endpoints $x_0$, $x_\infty$ are defined as follows. If $x_0$ (resp. $x_\infty$) is at the level $i$ and connected to the edge coming from the left, then the complex $P^{(*,0)}$ (resp. $P^{(*,\infty )}$) is
$$\ldots \rightarrow P_{i+2}\rightarrow P_{i+1} \rightarrow 0$$
with the differentials coming from the edges of the same color as the edge connecting to $x_0$ (resp. $x_\infty$).
If $x_0$ (resp. $x_\infty$) is at the level $i$ and connected to the edge coming from the right, then the complex $P^{(*,0)}$ (resp. $P^{(*,\infty )}$) is
$$\ldots \rightarrow P_{i+2}\rightarrow P_{i+1} \rightarrow P_i$$
with the differentials coming from the edges of the same color as the edge connecting to $x_0$ (resp. $x_\infty$).
The matrix entry of the differential from the summand $P_k$ corresponding to $e^{(k)}_l$ to the summand $P_i$ corresponding to $e^{(i)}_j$ is nonzero if the path from $e^{(i)}_j$ to $e^{(k)}_l$ goes directly to the right. The nonzero matrix entry is given by this path, multiplied by $-1$ if it is the lower path, and multiplied by the scalar $\lambda$ if the edge marked by $\lambda$ is on the path  from $e^{(i)}_j$ to $e^{(k)}_l$. Similarly for the matrix entries connecting the remaining parts of the resolution (corresponding to the end points $x_0$, $x_\infty$) to the rest.
\endproclaim

\demo{Proof} This is clear by definitions.
\qqed
\enddemo

Now we can finish the proof of the theorem A' . Consider the complex
$$\ldots\rightarrow Hom_{A(n)} (P^{(0)} , M(\beta  ;{\underline r};{\underline s}))\rightarrow Hom_{A(n)} (P^{(1)} , M(\beta  ;{\underline r};{\underline s})) \rightarrow 0.$$
We want to show that it is an isomorphism of vector spaces. We exhibit the bases of   spaces $Hom_{A(n)} (P^{(i)} , M(\beta  ;{\underline r};{\underline s}))$ similarly to the case of the band modules. Arrange basis vectors in $Hom_{A(n)} (P^{(i)} , M(\beta  ;{\underline r};{\underline s})))$ ordering the blocks in the way the extremal veritces are oredered in the band module, and within each block as we draw 
$\Gamma (\beta ;{\underline r};{\underline s})$. We construct a graph $EXT(M(\beta  ;{\underline r};{\underline s}))$ with the vertices of these spaces by connecting $e^{(i)}_u\otimes e^{(i)}_j$ to $e^{(k)}_v\otimes e^{(k)}_l$ if the matrix entry from the projectives cooresponding to the extremal vertices is nonzero, and the path given by matrix entry takes $e^{(i)}_u$ to $e^{(k)}_v$. 
Notice that the edges in the graph  $EXT(M(\beta  ;{\underline r};{\underline s}))$ consist of either upper or lower path in $\Gamma (\beta ; {\underline r}, {\underline s} )$. We refer to them as upper and lower edges in
$EXT(M(\beta  ;{\underline r};{\underline s}))$. The upper and lower edges in $EXT(M(\beta  ;{\underline r};{\underline s}))$ alternate, i.e. if a vertex has two edges adjacent to it, one is lower and the other is upper.
The main observation is that the graph $EXT(M(\beta  ;{\underline r};{\underline s}))$ has band and string components. Moreover, there are some distinguished band and string components with the vertices $e^{(i)}_j\otimes e^{(i)}_j$. Each block of vertices contains precisely one element of the distinguished component.

The remaining components are the even string components either occurring in degrees $0$ and $1$, or in degrees $>1$ -in the last case they are have two vertices. 
Indeed, we make the following observation about the nondistinguished  components occurring in degrees $0$ and $1$: either all their elements  occur above or below the  element from the distinguished component in their block.
The key point is that  if the orbit occurs above (resp. below) the band orbit in each block, then  its endpoints will be adjacent to the  upper (resp. lower) edge in $EXT(M(\beta  ;{\underline r};{\underline s}))$. Thus each string component has odd number of edges and thus even number of vertices.Thus the nondistinguished components  do not contribute to homology.

The distinguished components are band components except one string component.
This components is a string with odd number $2k+1$ verices with $k+1$ vertices in degree $0$, and $k$ vertices in degree $1$. Thus the only homology of the complex is 
$$Ext^0_{A(n)}(M(\beta  ;{\underline r};{\underline s}) , M(\beta  ;{\underline r};{\underline s}))=K.$$
This completes the proof. 
\qqed
\enddemo

\proclaim {Remark 5} Similar reasoning applies to arbitrary up-and-down module \break $M(\beta  ;{\underline r};{\underline s}))$, the only difference being we can have more than one distinguished string component. Thus higher Ext's will vanish also in such case.
\endproclaim

Let $C(\beta ;{\underline r}, {\underline s})$ be a string component. We denote by $N(\beta ;{\underline r}, {\underline s})$ the number of connected components of $\Gamma(\beta ;{\underline r}, {\underline s})$,
i.e. the number of summands. We have the following reduction which is a modfication of rules $a), b), c)$ given for band module components above.

We give an algorithm describing the  generic decomposition of the string components.

\proclaim {Theorem B'}  Assume that $C(\beta ;{\underline r};{\underline s})$ is a string component. Choose maximal $\beta_i$ and assume $\beta_{i-1}\ge \beta_{i+1}$. We can also assume that $\beta_{i+1}<\beta_i$. Define the number $t$ to be the difference between the horizontal change of red arrows and the level of blue arrows from the $(i-1)$-st to $i$-th level. More precisely, $t=\beta_i -\beta_{I-1}$. 
 We have two cases
\item{a)} $t>0$. Define the component $C(\beta ' ;{\underline r'};{\underline s'})$ by subtracting $t$ from $\beta_{i-1},\beta_i , r_i ,s_i $. Then $N(\beta ;{\underline r}; {\underline s})= N(\beta ' ;{\underline r'}; {\underline s'})$. Moreover, the summands from generic decompositions of $C(\beta ;{\underline r}; {\underline s})$ and  $C(\beta' ;{\underline r'};{\underline s'})$ differ only in places $i-1$ and $i$.
\item{b)} $t=0$. Define the component $C(\beta ' ;{\underline r'}; {\underline s'})$ by subtracting $1$ from $\beta_{i-1},\beta_i , r_i ,s_i $. Then $N(\beta ;{\underline r}; {\underline s})= N(\beta ' ;{\underline r'}; {\underline s'})+1$. Moreover, the generic decompositions of the components $C(\beta ;{\underline r}; {\underline s})$ and  $C(\beta ' ;{\underline r'};{\underline s'})$ are the same, with the generic decomposition of $C(\beta ;{\underline r}; {\underline s})$ containing additional summand of dimension vector 
$\epsilon_{i-1}+\epsilon_i$.
\endproclaim

\demo{Proof} One can prove the above reduction for the graphs $C(\beta ;{\underline r}; {\underline s})$ in similar way to the reduction for band components given above. We leave this to the reader.
\qqed
\enddemo

\head  \S 4.  The matching diagrams and associated toric varieties. \endhead

Let $D= (Z(n-1), \Theta )$ be a pair consisting of the set 
$$Z(n-1)=\lbrace x_1 ,\ldots ,x_{n-1}, y_1 ,\ldots ,y_{n-1})$$
 and the involution $\Theta : Z(n-1)\rightarrow Z(n-1)$ without fixed points.  We assume that $\Theta$ is {\it symmetric} , i.e. $\Theta (x_i )= x_j$ is equivalent to $\Theta (y_i )=y_j$, and $\Theta (x_i )= y_j$ is equivalent to $\Theta (y_i )=x_j$. Equivalently, symmetry means that $\Theta$ commutes with the involution $\eta$ defined by $\eta (x_i )=y_i$ and $\eta (y_i )=x_i$ for each $i$. We call the pair $D$ {\it a symmetric matching} (of two diagrams of type $A_n$).

We associate to the pair $D$ a toric variety as follows. Consider the polynomial ring
$$S=K[A_1 ,\ldots ,A_n , B_1 ,\ldots ,B_n ]$$
We define a subring $S(\Theta )$ as follows. The monomial $A_1^{u_1}\ldots A_n^{u_n}B_1^{v_1}\ldots B_n^{v_n}$ is in $S(\Theta)$ if the following equations, referred to as {\it the control equations},  are satisfied
$$u_i +u_{i+1}=v_j+v_{j+1}, \ \ \ \ u_j+u_{j+1}= v_i +v_{i+1}$$
if $\Theta (x_i )=y_j$, $\Theta (y_i )=x_j$ and
$$u_i +u_{i+1}=u_j+u_{j+1}, \ \ \ \  v_j+v_{j+1}= v_i +v_{i+1}$$
if $\Theta (x_i )=x_j$, $\Theta (y_i )=y_j$, for $1\le i,j\le n-1$.
We define an action of $(K^* )^{n-1}$ on $S$ as follows. Denote the orbits of $\Theta$ by $O_1 ,\ldots ,O_{n-1}$. Then the $j$-th copy of $K^*$ acts on $S$ as follows.

 If $O_j = \lbrace x_k ,y_l \rbrace$ then
$$t (A_1^{u_1}\ldots A_{n-1}^{u_{n-1}} B_1^{v_1}\ldots B_{n-1}^{v_{n-1}} )=  t^{u_k +u_{k+1}-v_l -v_{l+1}} 
A_1^{u_1}\ldots A_{n-1}^{u_{n-1}} B_1^{v_1}\ldots B_{n-1}^{v_{n-1}}.$$ 
 If $O_j = \lbrace x_k ,x_l \rbrace$ with $k<l$ then
$$t (A_1^{u_1}\ldots A_{n-1}^{u_{n-1}} B_1^{v_1}\ldots B_{n-1}^{v_{n-1}} )=  t^{u_k +u_{k+1}-u_l -u_{l+1}} 
A_1^{u_1}\ldots A_{n-1}^{u_{n-1}} B_1^{v_1}\ldots B_{n-1}^{v_{n-1}}.$$
Similarly if $O_j = \lbrace y_k ,y_l \rbrace$.

\proclaim {Proposition 7} $S(\Theta )=S^{(K^* )^{n-1}}$.
\endproclaim

\demo{Proof.} This follows directly from the definitions.
\qqed
\enddemo

\proclaim {Definition 4} The ring $S(\Theta)$ is called the matching ring associated to the matching $D$.
\endproclaim

\proclaim {Proposition 8} Consider the involution $\Theta (x_i )=y_i$. Then the subring
$S(\Theta )$ consists of all monomials $A_1^{u_1}\ldots B_n^{v_n}$ such that for each
$i$ we have $u_i +u_{i+1}=v_i +v_{i+1}$. In this case one can see easily that
$S(\Theta )$ is generated by $n+1$ elements $X_i = A_i B_i$ ($1\le i\le n$,
$Y_1= A_1 B_2 A_3 B_4\ldots$, $Y_2=B_1 A_2 B_3 A_4 \ldots$ with the relation
$X_1 \ldots X_{n} =Y_1Y_2$, so it is a hypersurface.
\endproclaim

Now we consider the involutions $\Theta$ with some special properties. These are important because below we will show that 
the involutions coming from the algebras $A(n)$ satisfy these properties.

\proclaim {Definition 5} Let $D = (Z(n-1), \Theta )$ be a symmetric matching.  
We say that numbers $i$ and $j$ are linked if $\Theta (x_i )= x_j$ or $\Theta (x_i )= y_j$.
The matching $D$ is even  if $\Theta (x_i )=y_j$ 
implies that $i-j$ is even and  $\Theta (x_i )=x_j$ implies that $i-j$ is odd.
The matching $D$ is unmixed if it cannot happen that there are pairs $(i,k)$ and $(j,l)$ of linked numbers that are ''mixed", i.e. such that $i<j<k<l$.
\endproclaim

\proclaim {Proposition 9} The ring $S(\Theta )$ is a complete intersection
for all symmetric unmixed even matchings $D= (Z(n-1), \Theta )$, for $n\le 6$. 
\endproclaim

\demo{Proof.} Let us start with a reduction lemma.

\proclaim {Lemma 2} Assume that the control equations include the equations 
$$u_{i-1} +u_i =u_i +u_{i+1},\ \ \ \  v_{i-1}+v_i =v_i +v_{i+1}.$$
for som $2\le i\le n-1$. 
 Then the ring $S(\Theta )$ is isomorphic to a polynomial ring in two variables over a matching ring 
$S(\Theta ')$ where we get $\Theta '$ by setting
$u'_j =u_j ,v'_j =v_j$ for $1\le j\le i-1$, and $u'_j =u_{j-2}, v'_j =v_{j-2}$ for $j\ge i+1$.
\endproclaim

\demo{Proof.} Notice that the control equation for $\Theta$ imply $u_{i-1}=u_{i+1}, v_{i-1}=v_{i+1}$, and that $A_i \in S(\Theta ), B_i\in S(\Theta )$.
The remainig equations after reindexing are the same as control equations for $S(\Theta' )$. Thus $S(\Theta )=S(\Theta ')[A_i ,B_i ]$ and $A_i ,B_i$ are algebraically independent over $S(\Theta ')$.
\qqed
\enddemo

This means that checking property of the ring $S(\Theta )$ for a matching $\Theta$ satisfying the assumption of the lemma reduces to a smaller case.
Thus we will call such matching {\it reducible}.

The proof of Proposition 9 is a tedious case by case check. 

For $n=2$ we have only one case
$\Theta (x_1 )=y_1$ so the results follows from the Proposition 8.
For $n=3$ there are two cases; the involution $\Theta (x_i )=y_i$ which is covered by Proposition 8 and the involution
$\Theta (x_1) =x_2$, $\Theta (y_1 )=y_2$. In this case the ring $S(\Theta )$ is generated by $X_1=A_2, X_2=B_2 ,Y_1=A_1 A_3 ,Y_2=B_1 B_3$ so it is a polynomial ring.

For $n=4$ we have two non-reducible cases:
\item{a)} $\Theta (x_i )=y_i$, which is a hypersurface by the Proposition 8,
\item{b)} $\Theta (x_1 )=y_3$, $\Theta (x_2 )=y_2$, $\Theta (x_3)=y_1$. 
The equations are
$$u_1 +u_2 =v_3 +v_4 , \ \ \ \ u_2 +u_3 =v_2 +v_3 , \ \ \ \ u_3 +u_4 =v_1 +v_2 .$$
The ring $S(\Theta )$ has the presentation
$$S(\Theta )=K[X_1 ,X_2 ,X_3 , X_4, Y_1 , Y_2]/(Y_1Y_2-X_1X_2X_3X_4)$$
where 
$$X_1=A_1 B_4, \ \ \ \ X_2=A_2 B_3, \ \ \ \ X_3=A_3 B_2, \ \ \ \  X_4=A_4 B_1,$$
$$ Y_1= A_1 A_3 B_1 B_3,  \ \ \ \ Y_2= A_2 A_4 B_2 B_4 .$$

For $n=5$ we have four non-reducible cases:
\item{a)} $\Theta (x_i )=y_i$, which is a hypersurface by the Proposition 8,
\item{b)} $\Theta (x_1 )=y_3$, $\Theta (x_2 )=y_2$, $\Theta (x_3 )=y_1$, $\Theta (x_4 ) =y_4 .$
$$u_1 +u_2 = v_3 +v_4 , \ \ \ \ u_2 +u_3 = v_2 +v_3 , \ \ \ \ u_3 +u_4 = v_1 +v_2 , \ \ \ \ u_4 +u_5 = v_4 +v_5 .$$
The ring $S(\Theta )$ has the presentation
$$S(\Theta )=K[X_1,X_2,X_3,Y_1,Y_2,Z_1,Z_2,Z_3]/(X_1X_2Y_1Y_2-X_3Z_1Z_3, X_1X_2Z_2-Z_1Z_3)$$
$$X_1= A_2 B_3 ,\ \ \ \ X_2=A_3 B_2 ,\ \ \ \  X_3=A_5 B_5, \ \ \ \ Y_1= A_1 A_5 B_4 , \ \ \ \ Y_2= A_4 B_1 B_5 ,$$
$$Z_1= A_1 A_3 B_1 B_3 , \ \ \ \ \ Z_2 =A_1A_4B_1B_4, \ \ \ \ Z_3=A_2 A_4 B_2 B_4.$$
\item{c)} $\Theta (x_1 )=y_1$, $\Theta (x_2 )=y_3$, $\Theta (x_3 )=y_3$, $\Theta (x_4 ) =y_2 .$
This case is symmetric to b).
\item{d)} $\Theta (x_1 )=x_4$, $\Theta (x_2 )=y_2$, $\Theta (x_3 )=y_3$, $\Theta (y_1 ) =y_4 .$
$$u_1 +u_2 = u_4 +u_5 , \ \ \ \ u_2 +u_3 = v_2 +v_3 , \ \ \ \ u_3 +u_4 = v_3 +v_4 , \ \ \ \ v_1 +v_2 , = v_4 +v_5 .$$
The ring $S(\Theta )$ has the presentation
$$S(\Theta )=K[X_1,X_2,X_3,Y_1,Y_2,Z_1,Z_2,Z_3]/(X_1X_2Z_3-Z_1Z_2, X_3Z_3-Y_1Y_2)$$
$$X_1=A_1 A_5 ,\ \ \ \ X_2= B_1 B_5 ,\ \ \ \ X_3=A_3 B_3 ,\ \ \ \  Y_1= A_2 A_4 B_3, \ \ \ \ Y_2= A_3 B_2 B_4,$$
$$Z_1= A_1 A_4 B_1 B_4 ,\ \ \ \ Z_2= A_2 A_5 B_2 B_5 ,\ \ \ \ Z_3= A_2 A_4 B_2 B_4 ,$$
so it is a complete intersection of codimension 2.

For $n=6$ we have eight non-reducible cases.

\item{a)} $\Theta (x_i )=y_i$, which is a hypersurface by the Proposition 8,
\item{b)} $\Theta (x_1 )=y_3$, $\Theta (x_2 )=y_2$, $\Theta (x_3 )=y_1$, $\Theta (x_4 ) =y_4 ,$ $\Theta (x_5 )= y_5$.
$$S(\Theta )=K[X_1,X_2,X_3,X_4,Y_1,Y_2,Y_3,Y_4,Y_5]/(Y_1Y_3-X_1X_2Y_2, Y_4Y_5-X_3X_4Y_2)$$
$$X_1=A_2B_3 ,\ \ \ \ X_2=A_3B_2, \ \ \ \ X_3 =A_5B_5, \ \ \ \ X_4 =A_6B_6,$$
$$ Y_1=A_1A_3B_1B_3 ,\ \ \ \ Y_2 =A_1A_4B_1B_4 ,\ \ \ \ Y_3 =A_2A_4B_2B_4 ,$$
$$Y_4 =A_1A_5B_4B_6, \ \ \ \ Y_5 =A_4A_6B_1B_5 .$$
\item{c)} $\Theta (x_1 )=y_1$, $\Theta (x_2 )=y_4$, $\Theta (x_3 )=y_3$, $\Theta (x_4 ) =y_2 ,$ $\Theta (x_5 )= y_5$.
$$S(\Theta )=K[X_1,X_2,X_3,X_4,Y_1,Y_2,Y_3,Y_4,Y_5]/(X_1X_4Y_5-Y_3Y_4,X_2X_3Y_5-Y_1Y_2)$$
$$X_1=A_1B_1,  \ \ \ \ X_2=A_3B_4, \ \ \ \ X_3=A_4B_3, \ \ \ \ X_4=A_6B_6, $$
$$Y_1=A_2A_4B_2B_4, \ \ \ \ Y_2 =A_3A_5B_3B_5, \ \ \ \ Y_3=A_1A_5B_2B_6,$$
$$ Y_4=A_2A_6B_1B_5, \ \ \ \ Y_5=A_2A_5B_2B_5.$$
\item{d)} $\Theta (x_1 )=y_1$, $\Theta (x_2 )=y_2$, $\Theta (x_3 )=y_5$, $\Theta (x_4 ) =y_4 ,$ $\Theta (x_5 )= y_3$.
This case is symmetric to $b)$.
\item{e)} $\Theta (x_1 )=y_5$, $\Theta (x_2 )=y_4$, $\Theta (x_3 )=y_3$, $\Theta (x_4 ) =y_2 ,$ $\Theta (x_5 )= y_1$.
This case is symmetric to $a)$ by exchanging $B_i$ with $B_{6-i}$.
\item{f)} $\Theta (x_1 )=y_5$, $\Theta (x_2 )=y_2$, $\Theta (x_3 )=y_3$, $\Theta (x_4 ) =y_4 ,$ $\Theta (x_5 )= y_1$.
This case is symmetric to $c)$ by exchanging $B_i$ with $B_{6-i}$.
\item{g)} $\Theta (x_1 )=x_4$, $\Theta (x_2 )=y_2$, $\Theta (x_3 )=y_3$, $\Theta (x_4 ) =x_1 ,$ $\Theta (x_5 )= y_5$.
$$S(\Theta)=$$
$$K[X_1,X_2,Y_1,Y_2,Y_3,Y_4,Z_1,Z_2,Z_3,Z_4]/(Y_1Y_4-X_2Z_2,Y_2Y_3-X_1Z_3,Z_1Z_4-Z_2Z_3)$$
$$X_1=A_3B_3, \ \ \ \ X_2=A_6B_6,$$
$$Y_1=A_1A_5B_6,\ \ \ \ Y_2=A_2A_4B_3,\ \ \ \ Y_3=A_3B_2B_4,\ \ \ \ Y_4=A_6B_1B_5,$$
$$Z_1=A_1A_4B_1B_4,\ \ \ \ Z_2=A_1A_5B_1B_5,\ \ \ \ Z_3=A_2A_4B_2B_4,\ \ \ \ Z_4=A_2A_5B_2B_5.$$
\item{h)} $\Theta (x_1 )=y_1$, $\Theta (x_2 )=x_5$, $\Theta (x_3 )=y_3$, $\Theta (x_4 ) =y_4 ,$ $\Theta (x_5 )= x_2$.
This case is symmetric to $g)$.
This concludes the proof of Proposition 9.
\qqed
\enddemo

Next we have the crucial example.

\proclaim {Example 8} Let $n=7$ and let
$\Theta (x_1 )=y_3$, $\Theta (x_2 )=y_2$, $\Theta (x_3 )=y_1$, $\Theta (x_4 ) =y_6,$ $\Theta (x_5 )= y_5$, $\Theta (x_6 )=y_4$.
The control equations are
$$u_1+u_2 =v_3 +v_4,$$
$$u_2+u_3 =v_2 +v_3,$$
$$u_3+u_4 =v_1 +v_2,$$
$$u_4+u_5 =v_6 +v_7,$$
$$u_5+u_6 =v_5 +v_6,$$
$$u_6+u_7 =v_4 +v_5.$$
The generators are:
$$X_1 =A_2B_3, X_2=A_3B_2 , X_3 = A_5B_6, X_4 =A_6B_5 ,$$
$$Y_1=A_1A_7B_4 ,\ \ \ \ Y_2=A_4B_1B_7 ,$$
$$Z_1 = A_1A_3B_1B_3 , \ \ \ \ Z_2=A_5A_7B_5B_7 ,$$
$$W_1=A_1A_4A_6B_1B_4B_6 ,\ \ \ \ W_2 =A_2A_4A_6B_2B_4B_6, \ \ \ \ W_3 =A_2A_4A_7B_2B_4B_7 .$$
The relations are the $4\times 4$ pfaffians of the skew symmetric  matrix
$$\left(\matrix 0&0&-W_1&Y_1Y_2&-Z_1\\
0&0&-W_2&W_3&-X_1X_2\\
W_1&W_2&0&0&-X_3X_4\\
-Y_1Y_2&-W_3&0&0&Z_2\\
Z_1&X_1X_2&X_3X_4&-Z_2&0
\endmatrix\right)$$
so it is a codimension 3 Gorenstein ideal with 5 generators, hence not a complete intersection.
\endproclaim

We finish this section with a conjecture about the generators and relations of the matching rings $S(\Theta )$.

\proclaim {Conjecture} Let $D= (Z(n-1), \Theta )$ be a symmetric matching. 
\item{a)} The ring $S(\Theta )$ is generated by the monomials $A_1^{u_1}\ldots A_n^{u_n}B_1^{v_1}\ldots B_n^{v_n}$ such that $u_i\le 1$ and $v_i\le 1$ for $i=1,\ldots ,n$.
\item{b)} The relations among the generators of the ring $S(\Theta )$ occur in multidegrees\break  $(u_1 ,\ldots ,u_n , v_1 ,\ldots ,v_n )$ with $u_i\le 2$ and $v_i\le 2$ for $i=1,\ldots ,n$.
\endproclaim

\head  \S 5.  The rings of semiinvariants for the algebras $A(n)$. \endhead

We consider the band module components $\beta = [m_{n},\ldots ,m_1 ]$ defined in the previous section.
We  investigate the rings of semiinvariants 
$$SI(A(n),\beta ):= K[R(A(n)),\beta )]^{SL(\beta )}.$$
The main result of this section is  the explicit description of the rings $SI(A(n),\beta )$ by generators 
and relations.  We will also show that  for low values of $n$ these rings are always complete intersections. 

The usual approach to semiinvariants is through Schofield semiinvariants $c^V$  corresponding to the indecomposable representations of projective dimension $\le 1$. These were proven to generate the rings of semiinvariants on irreducible components of representation spaces in \cite{DW3}. 

However because of the special nature of algebras $A(n)$ it is more convenient to describe the semiinvariants in terms of partitions. This description actually suggests the definition of the up-and-down modules and allows to construct them directly from the root systems. This approach should generalize to arbitrary string algebras. 

Let $\beta = [m_n ,\ldots ,m_1 ]$ be a band module component. The representation space $R(A(n),\beta )$ is a product of two varieties of complexes. Using Proposition 2 we get 

\proclaim {Corollary 1} The variety $Rep( A(n), \beta )$ of representations of
$ A(n)$ in dimenson $\beta$ has a distinguished irreducible component $R(A(n),\beta )$ containing pairs of exact complexes. 
This component has rational singularities.The coordinate ring $K[ R(A(n), \beta )]$ has the following decomposition
into Schur functors.
$$K[ R(A(n),\beta ) ]=\oplus_{\lambda_1 ,\ldots ,\lambda_n ,\mu_1 ,\ldots ,\mu_n} S(\lambda, \mu)$$
where
$$S(\lambda, \mu)=\eqalign{&S_{\lambda_n}F_n \otimes S_{(\lambda_{n-1},-\lambda_n )}F_{n-1}\otimes\ldots\otimes S_{(\lambda_1 ,-\lambda_2 )}F_1\otimes S_{(-\lambda_1 ) }F_0\otimes \cr
\otimes&S_{\mu_n}F_n \otimes S_{(\mu_{n-1},-\mu_n )}F_{n-1}\otimes\ldots\otimes S_{(\mu_1 ,-\mu_2 )}F_1\otimes S_{(-\mu_1 ) }F_0}$$
and we sum over the $2n$-tuples of partitions $(\lambda_n ,\ldots ,\lambda_1 ,\mu_n ,\ldots ,\mu_1 )$ with
$\lambda_i$ and $\mu_i$ having $m_i$ parts.
\endproclaim

Let us denote $\lambda_i = (\lambda_{i,1},\ldots ,\lambda_{i,m_{i-1}})$, $\mu_i = (\mu_{i,1},\ldots ,\mu_{i,m_{i-1}})$.
We also denote $\Lambda_i = (\lambda_i ,-\lambda_{i+1})$, $M_i = (\mu_i ,-\mu_{i+1})$ for $i=0,\ldots ,n$.
We denote by $(\Lambda ,M )$ the highest weight $(\Lambda_n ,\ldots ,\Lambda_1 ,M_n ,\ldots ,M_1 )$. 

It is clear from Littlewood-Richardson rule (comp. [MD], ch.1, [W], ch.2) that the summand corresponding
to the $2n$-tuple $(\lambda_n ,\ldots ,\lambda_1 ,\mu_n ,\ldots ,\mu_1 )$ contains
at most one $SL(\beta )$-invariant. This summand contains an invariant if
and only if for each $i=1,\ldots ,n$ and each $j=1,\ldots ,m_i$ we have 
$\Lambda_{i,j}= M_{i,m_i +1-j}$. It is enough to find a criterion when this happens and
then deduce the structure of the ring $SI(Q, \beta)$.

Let $\Delta$ denote the set of simple roots for the group $SL(A(n), \beta )$. We denote
by $\hat\sigma$ the involution of $\Delta$ corresponding to the symmetry of each of Dynkin components of
type $A_n$, i.e. changing the root $\epsilon_{i,j} -\epsilon_{i,j+1}$ into
$\epsilon_{i, \beta_i -j}-\epsilon_{i, \beta_i +1-j}$. 

Let us denote by $\hat\theta$ the involution of $\Delta$ related to the symmetry of highest weights
$\Lambda ,M$ ,i.e. $\theta$ interchanges the roots $\epsilon_{i,j} -\epsilon_{i,j+1}\equiv \epsilon_{i+1, \beta_i -j}-\epsilon_{i+1,\beta_i +1-j}$
for all $i$, $1\le i<j\le m_{i-1} -1$.

Consider the set $\Sigma =\Delta\coprod\Delta$. Two copies of $\Delta$ will be
denoted by $\Delta_0 $ and $\Delta_1$ and so will be corresponding simple roots.

We have two involutions of the set $\Sigma$, which we will denote by $\sigma$ and $\theta$. They are given by the formulas
${\sigma} (\alpha_s )=({\hat\sigma}(\alpha ))_{1-s}$ for $s=0,1$;
${\theta} (\alpha_s)=({\hat\theta}\alpha )_s$ for $s=0,1$.
We will see a close relationship of $\theta$ and $\sigma$ to semiinvariants.

\proclaim {Remark 6} It is interesting to note that the set $\Sigma$ with involutions $\sigma$ and $\theta$ is closely related to up-and-down modules, in fact it suggested their definition. 
We can construct from $\Sigma$ a double colored graph $\Gamma$ as follows. The vertices of $\Gamma$ are just $\Delta_0$, for even $i$, with the vertex $\alpha_1$ identified with $\theta (\alpha_1 )\in\Delta_0$. For odd $i$ the vertices of $\Gamma$ will be $\Delta_1$, with the roots $\alpha_0$ identified with $\theta (\alpha_0 )\in\Delta_1$. The red and blue arrows will be assigned via involutions $\theta\sigma$ or $\sigma\theta$ according to whether we go from $i+1$ to $i$ for $i$ even or odd.
In this way we get the graph $\Gamma (\hat\beta ,{\underline r}, {\underline s})$ where ${\hat\beta }= (\beta_n -1,\ldots,\beta_0 -1 )$,
${\underline r}={\underline s} =(m_{n-1}-1,\ldots ,m_0 -1 )$.
\endproclaim

For each pair of highest weights $(\Lambda ,M )$ we introduce the 
 function  
$$f(\Lambda, M): \Sigma\rightarrow N$$
defined by
$$f(\Lambda ,M)(\alpha_0 )= (\Lambda \alpha ),$$
$$f(\Lambda ,M)(\alpha_1 )= (M \alpha ).$$

We want to parametrize the pairs $(\Lambda ,M)$ that contain the semiinvariant.
We define an equivalence relation $\equiv$
on the set $\Sigma$ as the smallest equivalence relation containing the equivalences
$\alpha_0 \equiv \sigma (\alpha )_1$ and $\alpha_s \equiv \theta (\alpha )_s$ for
$s=0,1$.

\proclaim {Proposition 10} The summand corresponding to a weight $(\Lambda ,M)$ contains a semiinvariant if and
only if for $\alpha\equiv\beta$ we have $f(\Lambda ,M)(\alpha )=f(\Lambda ,M)(\beta )$.
\endproclaim

\demo{Proof} The proof follows from the definitions in a straightforward way.
\qqed
\enddemo

Let us consider the set of all weights $(\Lambda ,M)$ for which the corresponding summand $S(\Lambda , M)$ contains a semiinvariant. These pairs form
a semigroup $\Sigma (A(n), \beta )$ with respect to addition. 

\proclaim {Lemma 3} The ring of semi-invariants $SI(A(n), \beta )$ is isomorphic to a semigroup ring
$K[\Sigma (A(n),\beta )]$.
\endproclaim

\demo{Proof} For each $(\Lambda , M)\in \Sigma (A(n),\beta )$ we define a semi-invariant
$C(\Lambda, M)$ as follows.
The product of maximal exterior powers of $F_i$ embeds by the obvious diagonalizations into the tensor product (each top exterior power diagonalizing to the product of exterior powers corresponding to matching parts of $\Lambda$ and $M$. Let us call this element $C'(\Lambda , M)$.
$$\eqalign{&\bigwedge^{\lambda_n}F_n \otimes \bigwedge^{(\lambda_{n-1},-\lambda_n )}F_{n-1}\otimes\ldots\otimes \bigwedge^{(\lambda_1 ,-\lambda_2 )}F_1\otimes \bigwedge^{(-\lambda_1 ) }F_0\otimes \cr\otimes&\bigwedge^{\mu_n}F_n \otimes \bigwedge^{(\mu_{n-1},-\mu_n )}F_{n-1}\otimes\ldots\otimes \bigwedge^{(\mu_1 ,-\mu_2 )}F_1\otimes \bigwedge^{(-\mu_1 ) }F_0}$$
covering the summand $S(\Lambda, M)$. 
However the tensor product of the exterior powers maps directly to the coordinate ring of our component via the map $m(\Lambda ,M)$ given by products of minors on the arrows of our quiver.
We define the semi-invariant $C(\Lambda, M)$ to be $m(\Lambda, M)C'(\Lambda , M)$. To show that these semi-invariants are linearly independent it is enough to show that the coset of $C'(\Lambda , M)$
is non-zero in
$$\eqalign{&S_{\lambda_n}F_n \otimes S_{(\lambda_{n-1},-\lambda_n )}F_{n-1}\otimes\ldots\otimes S_{(\lambda_1 ,-\lambda_2 )}F_1\otimes S_{(-\lambda_1 ) }F_0\otimes \cr
\otimes&S_{\mu_n}F_n \otimes S_{(\mu_{n-1},-\mu_n )}F_{n-1}\otimes\ldots\otimes S_{(\mu_1 ,-\mu_2 )}F_1\otimes S_{(-\mu_1 ) }F_0}$$
as this is a factor of a good filtration on the coordinate ring $K[R(A(n)),\beta )]$. Let $S(\Lambda ,M)_0$ be the span of $x\otimes y$ where $x$ is the highest weight vector in 
$S_{\lambda_n}F_n \otimes S_{(\lambda_{n-1},-\lambda_n )}F_{n-1}\otimes\ldots\otimes S_{(\lambda_1 ,-\lambda_2 )}F_1\otimes S_{(-\lambda_1 ) }F_0$
and $y$ is the lowest weight vector in
$S_{\mu_n}F_n \otimes S_{(\mu_{n-1},-\mu_n )}F_{n-1}\otimes\ldots\otimes S_{(\mu_1 ,-\mu_2 )}F_1\otimes S_{(-\mu_1 ) }F_0$.
The semi-invariant $C'(\Lambda ,M)$ has one non-zero summand in
$S(\Lambda ,M)_0$. We need to see that when straightening to write $C'(\Lambda ,M)$ as a combination of standard tableux this summand does not cancel out with another one.
 Howeve in straightening the double tableaux the weights of both components are preserved, and both the highest weight space and lowest weight space are one dimensional
 so the summand $x\otimes y$ does not cancel out and $C'(\Lambda ,M)$ is non-zero.
\qqed
\enddemo

In view of the lemma in the rest of the paper we identify  $SI(A(n),\beta )$ with the semigroup ring. Note that associating
to each pair $(\Lambda ,M)$ the function $f(\Lambda ,M)$ is a homomorphism of monoids.

We will first analyse the semigroup of functions $f(\Lambda ,M)$. This semigroup 
is  isomorphic to $N^x$ where $x$ is the number of all equivalence classes of $\equiv$. 

Let us investigate the equivalence classes of the relation $\equiv$.

We define the subset of critical elements of $Z(\underline m )$ to be the subset of roots of type
$(\epsilon_{i, m_{i-1}}-\epsilon_{i,m_{i-1}+1})_s$
for $i=1,\ldots ,n-1$, $s=0,1$. There are precisely $2(n-1)$ critical elements.

The equivalence class of $\equiv$
is called critical if it contains a critical element.

\proclaim {Proposition 11} There are exactly $n-1$ critical orbits, each of which contains exactly
two critical elements. 
\endproclaim

\demo {Proof} Let us look at the orbit of an element $\alpha_s \in\Sigma$.
From noncritical element $\alpha_s$ we can produce its $\theta$-neighbor 
 $\theta (\alpha )_s$ and its $\sigma$-neighbor $\sigma (\alpha )_{1-s}$.
We can get every element in the orbit by iterating these procedures, i.e taking
consecutive neighbors (alternating $\sigma$ and $\theta$).

Every element has a unique $\theta$-neighbor. The critical elements are the only ones which do not 
have a $\sigma$-neighbor, all other elements have a unive $\sigma$-neighbor.

Now the Proposititon is clear. 
\qqed
\enddemo

Let $O\subset\Sigma$ be a noncritical orbit of $\equiv$. Let's describe the pairs $(\Lambda ,M)$ for
which $f(\Lambda ,M)$ is the characteristic function of the orbit $O$.

\proclaim {Proposition 12}
 Let $O$ be a noncritical orbit. Then there is exactly one pair
$(\Lambda ,M)$ for which $f(\Lambda ,M)$ is the characteristic function of $O$.
\endproclaim

\demo {Proof} Let $O$ be a noncritical orbit. Let $(\Lambda ,M )$ be such
a pair of weights that $f(\Lambda ,M)$ is the characteristic function
of the orbit $O$. Rewriting the conditions in terms of $2n$-tuple of partitions
$(\lambda_n ,\ldots ,\lambda_1 ,\mu_n ,\ldots ,\mu_1 )$ we notice that since the values
of $f(\Lambda ,M)$ on critical elements are $0$, last part of each of the partitions is equal to $0$.
Since the values of $f(\Lambda ,M)$ on all elements of $\Sigma$ determine the difference of
all consecutive parts of each partition, we see that all partitions are determine uniquely.
\qqed
\enddemo

Let us denote $u_i = \Lambda_{i,m_i}$ and $v_i = M_{i,m_i}$. The values of the function $f(\Lambda ,M)$
determine all values $x_i =u_i +u_{i+1}$ and $y_i =v_i + v_{i+1}$. Each pair of numbers corresponding
to the points in the same critical orbit will be equal.

Let us associate to the data $\underline m =\lbrace m_{n-1} ,\ldots ,m_0 \rbrace$ the diagram
$D(\underline m )$ as follows. $D(\underline m ) = (Z,\Theta (\underline m ) )$ 
consists of a $2(n-1)$-tuple $Z = \lbrace x_1 ,\ldots ,y_{n-1}\rbrace$ of critical points
  together with involution $\Theta (\underline m ) :Z\rightarrow Z$
which sends a point $z\in Z$ to a point corresponding to the other end of the critical orbit
containing the point $z$. The involution $\Theta (\underline m )$ has the following properties.

\proclaim {Proposition 13} \item {a)} $\Theta (\underline m )$ is symmetric,i.e. if 
$\Theta (\underline m )(x_i )=x_j$ then
$\Theta (\underline m )(y_i )=y_j$, and if $\Theta (\underline m )(x_i )=y_j$ then  
$\Theta (\underline m )(y_i )=x_j$,
\item {b)} $\Theta (\underline m )$ is even, i.e. if $\Theta (\underline m )(x_i )=y_j$ 
then $i-j$ is even and 
if $\Theta (\underline m )(x_i )=x_j$ then $i-j$ is odd.
\endproclaim

\demo {Proof} The proof follows easily from the definition of the critical orbits.
\qqed
\enddemo

In the previous section we associated to the symmetric  involution $\Theta (\underline m )$ the subring
$S(\Theta (\underline m ))$.

This allows us to state the main result of the paper.

\proclaim {Theorem C} Let $\underline m =(m_{n-1} ,\ldots ,m_0 )$ be an $n$-tuple of natural numbers.
Let $\beta :=\beta ([m_{n},\ldots ,m_1])$ be the corresponding band module component for the algebra $A(n)$. 
Let $Z$ be the set of critical points
and let $\Theta (\underline m )$ be the involution of $Z$ determined by
the critical orbits. The ring of semiinvariants $SI(A(n), \beta)$
is isomorphic to the polynomial ring over the ring $S(\Theta (\underline m ))$.
Thus all the rings of semiinvariants are the polynomial rings over the rings of invariants of tori.
\endproclaim

\demo {Proof} By Lemma 3 the ring of semiinvariants $SI(A(n), \beta )$ is isomorphic
to the semigroup of pairs $(\Lambda ,M)$ containing a semiinvariant. Let us consider the function
$f(\Lambda ,M)$ introduced above. It's value on noncritical orbits can be set in an arbitrary
way using the polynomial generators described in Proposition 11. In terms of partitions
$(\Lambda ,M)$ it determines certain jumps in the partitions.

In order to determine the remaining parts of partitions from $(\Lambda ,M)$ notice
that it is enough to determine the last parts $u_i = \Lambda_{i,m_{i-1}}$ and $v_j = M_{j,m_{j-1}}$.
This will determine the values on the endpoints of critical orbits, and therefore
the remaining jumps in all partitions. The only constraint is that the values on the endpoints of
an orbit should be the same. This translates to the condition that the numbers
$u_i +u_{i+1}$ (or $v_j +v_{j+1}$) corresponding to the endpoints of a critical orbit
should be the same. This gives a homomorphism from the ring $S(\Theta (\underline m ))$ to
$SI(A(n), \beta )$ sending the monomial $(*)$ to the pair $(\Lambda ,M)$ such that
$u_i = \Lambda_{i,m_{i-1}}$ and $v_j = M_{j,m_{j-1}}$ and such that the jumps from noncritical orbits are $0$.
Then $SI(A(n), \beta)$ is clearly a polynomial ring over $S(\Theta (\underline m ))$.
\qqed
\enddemo
  
The involutions $\Theta (\underline m )$ one obtains have one more property. In order to establish it we investigate
 how the reduction given in Theorem B works for the rings of semiinvariants.

\proclaim {Theorem D} Let us consider a reduction step in the procedure given in Theorem B. Let $\beta :=\beta ([m_{n},\ldots ,m_1 ])$.
\item {a)} Assume that $m_{i+1} =m_i$ and consider the band module component 
$$\beta' :=\beta ([m_{n-1},\ldots ,m_{i+1}, m_{i-1},\ldots ,m_0]).$$
Then the involution $\Theta :=\Theta ([m_{n},\ldots ,m_1 ])$ is obtained from the involution 
$\Theta' :=\Theta ([m_{n-1},\ldots ,m_{i+1}, m_{i-1},\ldots ,m_0 ]),$
 by setting
$\Theta (x_{i+1})=y_{i+1},$
and the other values of the involution $\Theta$ are determined by the values of the involution $\Theta' $ by identifying the vertices $x_1 ,\ldots ,x_{n-2}$ (resp. the vertices $y_1 ,\ldots ,y_{n-2}$) of
the set $D([m_{n-1},\ldots ,m_{i+1}, m_{i-1},\ldots ,m_0 ])$ with the vertices $x_1 ,\ldots ,x_i ,x_{i+2},\ldots ,x_{n-1}$ (resp. the vertices $y_1 ,\ldots ,y_i ,y_{i+2},\ldots ,y_{n-1}$) of $D([m_{n},\ldots ,m_1 ])$,
\item{b)} Assume that $m_i >max (m_{i-1}, m_{i+1})$. Then the involution $\Theta ([m_{n},\ldots ,m_1])$ is the same as  $\Theta ([m_{n-1},\ldots ,m_{i+1}, m_i -|m_{i-1}-m_{i+1}|, m_{i-1},\ldots ,m_0])$, 
\item{c)} Assume that $m_i >m_{i-1}=m_{i+1}$. Then  the involution  $\Theta ([m_{n},\ldots ,m_1])$ satisfies $\Theta ([m_{n},\ldots ,m_1])(x_i)= x_{i+1}$,
$\Theta ([m_{n},\ldots ,m_1])(y_i)= y_{i+1}$. We get the involution $\Theta ([m_{n-1},\ldots ,m_{i+1}, m_{i+1}, m_{i-1},\ldots ,m_0])$ by making $x_i$ go to $y_i$, $x_{i+1}$
going to $y_{i+1}$ and keeping other values fixed.
\endproclaim

\demo{Proof} This follows from considering the up-and-down diagrams $\Gamma ([m_{n},\ldots ,m_1 ])$ and observing that if among edges from level $i+1$ to level $i$ we delete the lowest arrow from the upper part (i.e. red arrow for $n-i+1$ even and blue arrow for $n-i+1$ odd), the end points of remaining strings (except the lowest string of dimension vector $(1,1,\ldots ,1$ will give us the structure of $\Theta ([m_{n},\ldots ,m_1 ])$.
\qqed
\enddemo

The important corollaries are as follows. We say that numbers $i$ and $j$ are linked if $\Theta (x_i )= x_j$ (if $j-i$ is odd) or $\Theta (x_i )= y_j$ (if $j-i$ is even).

\proclaim {Corollary 2} 
\item{a)} The involutions $\Theta ([m_{n},\ldots ,m_1 ])$ have the following "unmixedness" property.
It cannot happen that there are pairs $(i,k)$ and $(j,l)$ of linked numbers that are ''mixed", i.e. such that $i<j<k<l$.
\item{b)} If a band module component is a Schur component (i.e. the corresponding up-and-down graph is connected) then the involution
$\Theta ([m_{n},\ldots ,m_1 ])$ is just given by the rule $\Theta ([m_{n},\ldots ,m_1 ])(x_i )=y_i$. Thus the corresponding semi-invariant ring $SI(A(n),\beta ([m_{n},\ldots ,m_1 ])$ is a hypersurface.
\endproclaim

\demo{Proof} Just follow the reduction procedure and investigate its effect on $\Theta$.
\qqed
\enddemo

\proclaim {Proposition 14} The ring $S(\Theta )$ is a complete intersection
for all symmetric unmixed even involutions $\Theta$, for $n\le 6$. Therefore the ring
$SI(A(n), \beta )$ is a complete intersection for all $n\le 6$. For $n=7$ the ring $S(\Theta (\underline m ))$
does not have to be a complete intersection, and so the same is true for the rings of semi-invariants
$SI(A(n),\beta )$ of band module components.
\endproclaim

\demo{Proof} 
This follows from Proposition 9 and Example 8. The only thing one needs to show is that the
matching $\Theta$ from Example 8 occurs as a matching for some band component.
In fact it is easy to check that this occurs for the dimension vector $\beta ([m_1,\ldots ,m_7])$ for
$ [m_1,\ldots ,m_7]= [1,2,2,1,2,3,1]$.

\qqed
\enddemo


\enddocument